\documentclass[12pt]{article}
\usepackage[utf8]{inputenc}

\usepackage[left=2cm, right=2cm, top=2.50cm, bottom=2.50cm]{geometry}

\usepackage{amssymb}
\usepackage[all,arc]{xy}
\usepackage{mathrsfs}
\usepackage{extarrows}
\usepackage{xcolor}
\usepackage{amsmath}
\usepackage{amsfonts,epsfig}
\usepackage{amssymb,bm}
\usepackage{graphicx}
\usepackage{subfig}
\usepackage{pgf,pgfarrows,pgfnodes,pgfautomata,pgfheaps}
\usepackage{appendix}

\usepackage{mathtools}
\usepackage{tabularx}  
\usepackage{booktabs}
\usepackage{multirow}

\usepackage{array,makecell}

\numberwithin{equation}{section}
\newenvironment {Proof} {\noindent {\em Proof.}}

\usepackage{amsmath,amsthm,amsfonts,amssymb,amscd,color}
\usepackage{hyperref}
\hypersetup{hypertex=true,
colorlinks=true,
linkcolor=blue,
anchorcolor=blue,
citecolor=blue}

\allowdisplaybreaks[4]
\newtheorem{theo}{Theorem}[section]
\newtheorem{lem}[theo]{Lemma}
\newtheorem{cor}[theo]{Corollary}
\newtheorem{prop}[theo]{Proposition}

\newtheorem{remark}[theo]{Remark}

\newtheorem{que}[theo]{Question}

\begin{document}

\newcommand{\x}{\textbf{x}}
\newcommand{\y}{\textbf{y}}
\newcommand{\bu}{\textbf{b}_u}
\newcommand{\bv}{\textbf{b}_v}
\newcommand{\bubu}{\left \langle  \textbf{b}_u ,\textbf{b}_u \right \rangle}
\newcommand{\bubv}{\left \langle  \textbf{b}_u ,\textbf{b}_v \right \rangle}
\newcommand{\buj}{\left \langle  \textbf{b}_u ,\textbf{j}\ \right \rangle}

\title{Regular graphs with a complete bipartite graph as a star complement\thanks{This work is supported by the National Natural Science Foundation of China (Grant No. 11971180), the Guangdong Provincial Natural Science Foundation (Grant No. 2019A1515012052).}}
\author{Xiaona Fang$^a$, Lihua You$^a$\thanks{Corresponding author: ylhua@scnu.edu.cn}, Rangwei Wu$^a$, Yufei Huang$^b$\\
{\small $^a$School of Mathematical Sciences, South China Normal University,}\\ {\small Guangzhou, 510631, P. R. China}\\
{\small $^b$Department of Mathematics Teaching, Guangzhou Civil Aviation College,}\\ {\small Guangzhou, 510403, P. R. China}\\
}
\date{}
\maketitle

{\bf Abstract}: Let $G$ be a graph of order $n$ and $\mu$ be an adjacency eigenvalue of $G$ with multiplicity $k\geq 1$. A star complement $H$ for $\mu$ in $G$ is an induced subgraph of $G$ of order $n-k$ with no eigenvalue $\mu$, and the vertex subset $X=V(G-H)$ is called a star set for $\mu$ in $G$. The study of star complements and star sets provides a strong link between graph structure and linear algebra. In this paper, we study the regular graphs with $K_{t,s}\ (s\geq t\geq 2)$ as a star complement for an eigenvalue $\mu$, especially, characterize the case of $t=3$ completely, obtain some properties when $t=s$, and propose some problems for further study.



\baselineskip=0.30in

\section{Introduction}

\hspace{1.5em}
Let $G$ be a simple graph with vertex set $V(G)=\left \{ 1,2,\dots,n \right \}=\left [n\right ]$ and edge set $E(G)$. The adjacency matrix of $G$ is an $n\times n$ matrix $A(G)=(a_{ij})$, where $a_{ij}=1$ if vertex $i$ is adjacency to vertex $j$, and $0$ otherwise. We use the notation $i\sim j$ ($i\nsim j$) to indicate that $i,j$ are adjacent (not-adjacent) and the notation $d_G (i)$ to indicate the degree of vertex $i$ in $G$. The adjacency eigenvalues of $G$ are just the eigenvalues of $A(G)$. For more details on graph spectra, see \cite{1}.

Let $\mu$ be an eigenvalue of $G$ with multiplicity $k$. A {\em star set} for $\mu$ in $G$ is a subset $X$ of $V(G)$ such that $\left | X \right |=k$ and $\mu$ is not an eigenvalue of $G-X$, where $G-X$ is the subgraph of $G$ induced by $\overline{X}=V(G)\setminus X $. In this situation $H=G-X$ is called a {\em star complement} corresponding to $\mu$. Star sets and star complements exist for any eigenvalue of a graph, and they need not to be unique. The basic properties of star sets are established in Chapter $7$ of \cite{2}.

There is another equivalent geometric definition for star sets and star complements. Let $G$ be a graph with vertex set $V(G)=\left\{1,\dots,n \right\}$ and adjacency matrix $A$. Let $\left\{e_1,\dots,e_n\right\}$ be the standard orthonormal basis of $\mathbb{R}^n$, $\mu$ be an eigenvalue of $G$, and $P$ be the matrix which represents the orthogonal projection of $\mathbb{R}^n$ onto the eigenspace $\mathcal{E} (\mu)=\left\{ x \in \mathbb{R}^n: A(G)x=\mu x \ \right\}$ of $A$ with respect to $\left\{e_1,\dots,e_n\right\}$. Since $\mathcal{E} (\mu)$ is spanned by the vectors $Pe_j(j= 1,\dots, n)$, there exists $X \subseteq V (G)$ such that the vectors $Pe_j(j \in X)$ form a basis for $\mathcal{E} (\mu)$. Such a subset $X$ of $V(G)$ is called a {\em star set} for $\mu$ in $G$. In this situation $H=G-X$ is called a {\em star complement} for $\mu$.

For any graph $G$ of order $n$ with distinct eigenvalues $\lambda_1,\dots,\lambda_m$, there exists a partition $V(G)=X_1 \bigcup \cdots \bigcup X_m$ such that $X_i$ is a star set for eigenvalue $\lambda_i \ (i = 1,\dots,m)$. Such a partition is called a {\em star partition} of $G$.	For any graph $G$, there exists at least one star partition (\cite{5}). Each star partition determines a basis for $\mathbb{R}^n$ consisting of eigenvectors of an adjacency matrix. It provides a strong link between graph structure and linear algebra.


There are a lot of literatures about using star complements to construct and characterize certain graphs(\cite{krrr,bell2,bell,3,fang2021,6,row2,8,9,11,unicyclic,yuan3,jianfeng,yuan2,yuan}), especially, regular graphs with a prescribed graph such as $K_{1,s}$, $K_1\triangledown hK_q$, $K_{2,5}$, $K_{2,s}$, $K_{1,1,t}$, $K_{1,1,1,t}$, $\overline{sK_1 \cup K_t}$, $P_t \ (\mu=1)$, $K_{r,r,r}\ (\mu=1)$ or $K_{r,s}+tK_1\ (\mu=1)$ as a star complement were well studied in the literature. 
Motivated by the above research, in this paper, we introduce the fundamental properties of the theory of star complements in Section \ref{prelimi}, study the regular graphs with the bipartite graph $K_{t,s}\ (s\geq t\geq 1)$ as a star complement for an eigenvalue $\mu$ in Section \ref{kts}, completely characterize the regular graphs with $K_{3,s}\ (s\geq 3)$ as a star complement for an eigenvalue $\mu$ in Section \ref{regular}, study some properties of $K_{s,s}$ in Section \ref{kss}, and propose some problems for further research. 

\section{Preliminaries}\label{prelimi}

\hspace{1.5em}In this section, we introduce some results of star sets and star complements that will be required in the sequel.
The following fundamental result combines the Reconstruction Theorem (\cite[Theorem 7.4.1]{2}) with its converse (\cite[Theorem 7.4.4]{2}).

\begin{theo}{\rm(\cite{2})}\label{thm21}
	Let $\mu$ be an eigenvalue of $G$ with multiplicity $k$, $X$ be a set of vertices in the graph $G$.  Suppose that $G$ has adjacency matrix
	$$\left (\begin{matrix}
	A_X & B^T \\
	B & C
	\end{matrix}\right ),$$
	where $A_X$ is the adjacency matrix of the subgraph induced by $X$. Then $X$ is a star set for $\mu$ in $G$
	if and only if $\mu$ is not an eigenvalue of $C$ and
	\begin{equation}\label{main}
	\mu I-A_X=B^T(\mu I-C)^{-1}B.
	\end{equation}
	In this situation, $\mathcal{E}(\mu)$ consists of the vectors
	\begin{equation}\label{vector}
	\left ( \begin{matrix}
	\textbf{x} \\
	(\mu I-C)^{-1}B \textbf{x}
	\end{matrix}
	\right),
	\end{equation}
	where $\textbf{x} \in \mathbb{R}^k$.
\end{theo}

Note that if $X$ is a star set for $\mu$, then the corresponding star complement $H(=G-X)$ has adjacency matrix $C$,
and (\ref{main}) tells us  that $G$ can be determined by $\mu$, $H$ and the $H$-neighbourhood of vertices in $X$,
where the $H$-neighbourhood of the vertex $u\in X$, denoted by $N_H(u)$, is defined as $N_H(u)=\left\{v\mid v\sim u, v\in V(H)\right\}$.

It is usually convenient to apply (\ref{main})  in the form
\begin{equation*}
m(\mu)(\mu I-A_X)=B^Tm(\mu)(\mu I-C)^{-1}B,
\end{equation*}
where $m(x)$ is the minimal polynomial of $C$. This is because $m(\mu)(\mu I-C)^{-1}$ is given explicitly as follows.

\begin{prop}{\rm (\cite{3}, Proposition 0.2)}\label{prop26}
	Let $C$ be a square matrix with minimal polynomial
	$$m(x)=x^{d+1}+c_d x^d+c_{d-1} x^{d-1}+ \cdots +c_1 x+c_0.$$
	If $\mu$ is not an eigenvalue of $C$, then
	$$m(\mu)(\mu I-C)^{-1}=a_d C^d+a_{d-1} C^{d-1}+\cdots +a_1 C+a_0 I,$$
	where $a_d=1$ and for $0<i\leq d$,
	$a_{d-i}=\mu ^i+c_d \mu^{i-1}+c_{d-1} \mu^{i-2}+\cdots+c_{d-i+1}.$
\end{prop}

In order to find all the graphs with a prescribed star complement $H$ for $\mu$, we need to find all solution $A_X$, $B$ for given $\mu$ and $C$.
For any $\x, \y \in \mathbb{R}^q$, where $q=\left|V(H) \right|$, let
\begin{equation}\label{eq0022}
\left \langle  \x ,\y \right \rangle =\x ^T(\mu I-C)^{-1} \y.
\end{equation}

Let $\bu$ be the column of $B$ for any $u \in X$. By Theorem \ref{thm21}, we have

\begin{cor}{\rm (\cite{5}, Corollary 5.1.8 )}\label{cor22}
	Suppose that $\mu$ is not an eigenvalue of the graph $H$, where $\left|V(H)\right|=q$. There exists a graph $G$ with a star set $X=\left\{u_1,u_2,\dots,u_k\right\}$ for $\mu$ such that $G-X=H$ if and only if there exist $(0,1)$-vectors $\textbf{b}_{u_1},\textbf{b}_{u_2},\dots,\textbf{b}_{u_k}$ in $\mathbb{R}^q$ such that\\
	{\rm(1)} $\bubu =\mu$ for all $u \in X$, and \\
	{\rm (2)} $\bubv =\left\{\begin{matrix}
	-1, & u \sim v \\
	0, & u \nsim v
	\end{matrix}\right.$ for all pairs $u,v$ in $X$.
\end{cor}

In view of the two equations in the above corollary, we have

\begin{lem}{\rm (\cite{2})}\label{lem23}
	Let $X$ be a star set for $\mu$ in $G$, and  $H=G-X$.\\
	{\rm(1)} If $\mu \neq 0$, then $V(H)$ is a dominating set for $G$, that is, the $H$-neighbourhood of any vertex in $X$ are non-empty;\\
	{\rm(2)} If $\mu \notin \left\{ -1,0 \right\}$, then $V(H)$ is a location-dominating set for $G$, that is, the $H$-neighbourhood of distinct vertices in $X$ are distinct and non-empty.
\end{lem}

It follows from $(2)$ of Lemma \ref{lem23} that there are only finitely regular graphs with a prescribed star complement for $\mu \notin \left\{-1,0\right\}$. If $\mu =0$ and $X$ has distinct vertices $u$ and $v$ with the same neighbourhood in $G$, then $u$ and $v$ are called \emph{duplicate vertices}. If $\mu =-1$ and $X$ has distinct vertices $u$ and $v$ with the same closed neighbourhood in $G$, then $u$ and $v$ are called \emph{co-duplicate vertices} (see \cite{7}).

Recall that if the eigenspace $\mathcal{E}(\mu)$ is orthogonal to the all-1 vector $\textbf{j}$ then $\mu$ is called a \emph{non-main} eigenvalue. From (\ref{vector}), we have the following result.
\begin{lem}{\rm (\cite{3}, Proposition 0.3)}
The eigenvalue $\mu$ is a non-main eigenvalue if and only if
\begin{equation}\label{eq23}
    \buj =-1 \ \mbox{ for all } u \in X,
\end{equation}
where $\textbf{j}$ is the all-1 vector.
\end{lem}

\begin{lem}{\rm (\cite{5}, Corollary 3.9.12)}\label{lem25}
In an $r$-regular graph, all eigenvalues other than $r$ are non-main.
\end{lem}

In the rest of this paper, we let $H\cong K_{t,s}\ (s\geq t\geq 1)$, $(V, W)$ be a bipartition of the graph $K_{t,s}$ with $V=\left\{ v_1,v_2,\dots, v_t \right\}$, $W=\left\{ w_1,w_2,\dots,w_s \right\}$. We say that a vertex $u \in X$ is of type $(a,b)$ if it has $a$ neighbours in $V$ and $b$ neighbours in $W$. Clearly $(a,b)\neq (0,0)$ and $0\leq a\leq t$, $0\leq b\leq s$.

Let $C$ be the adjacency matrix of $H$, then $C$ has minimal polynomial $m(x)=x(x^2-ts).$
Since $\mu$ is not an eigenvalue of $C$, we have $\mu \neq 0$ and $\mu^2 \neq ts$. From Proposition \ref{prop26}, we have
\begin{equation}\label{eq24}
    m(\mu)(\mu I-C)^{-1}=C^2+\mu C+(\mu^2-ts)I.
\end{equation}
If $\mu$ is a non-main eigenvalue of $G$, then by (\ref{eq23}) and (\ref{eq24}) we have
\begin{equation}\label{eq25}
  \mu ^2(a+b)+\mu (as+tb)=-\mu (\mu^2-ts).
\end{equation}

Using (\ref{eq24}) to compute $\bubu =\mu$, we obtain the following equation
\begin{equation}\label{eq26}
    (\mu ^2-ts)(a+b)+a^2 s+tb^2+2ab \mu=\mu^2(\mu^2-ts).
\end{equation}

Let $u,v$ be distinct vertices in $X$ of type $(a,b)$, $(c,d)$, respectively. Let $\rho_{uv}=\left |N_H(u)\cap N_H(v)\right|$, $a_{uv}=1$ or $0$ according as $u\sim v$ or $ u\nsim v$. Using (\ref{eq24}) to compute $\bubv =-a_{uv}$, we have
\begin{equation}\label{eq27}
    (\mu^2-ts)\rho_{uv}+acs+bdt+\mu (ad+bc)=-\mu (\mu^2-ts)a_{uv}.
\end{equation}

\section{Regular graphs with $K_{t,s}$ as a star complement}\label{kts}

\hspace{1.5em}An $r$-regular graph $G$ with $n$ vertices is said to be strongly regular with parameters $(n,r,e,f)$ if every two adjacent vertices in $G$ have $e$ common neighbours and every two non-adjacent vertices have $f$ common neighbours. For example, Petersen graph is strongly regular with parameters $(10,3,0,1)$.

For the regular graphs with the complete bipartite graph $K_{t,s}$ as a star complement, the case of $t=1$ was solved by Rowlinson and Tayfeh-Rezaie in 2010 (\cite{9}), the case of $t=2,\ s=5$ was solved by Rowlinson and Jackson in 1999 (\cite{8}), the case of $t=2,\ s\neq 5$ was solved by Yuan, Zhao, Liu and Chen in 2018 (\cite{yuan}), and the conclusions are listed below.

\begin{theo}\rm{(\cite{9})}\label{thm41}	
	If the $r$-regular graph $G$ has $K_{1, s}\ (s>1)$ as a star complement for an eigenvalue $\mu$, then one of the following holds:
	\\(1) $\mu=\pm 2, r=s=2$ and $G\cong K_{2,2}$;
	\\(2) $\mu=\frac{1}{2}(-1 \pm \sqrt{5}), r=s=2$ and $G$ is a 5-cycle;
	\\(3) $\mu \in \mathbb{N}_+, r=s$ and $G$ is strongly regular with parameters $\left(\left(\mu^{2}+3 \mu\right)^{2}, \mu\left(\mu^{2}+3 \mu+1\right), 0, \mu(\mu+1)\right)$.
\end{theo}

\begin{theo}\rm{(\cite{8,yuan})}\label{thm42}
	Let $s \geq 2$. If the $r$-regular graph $G$ has $K_{2, s}$ as a star complement for an eigenvalue $\mu$, then one of the following holds:
	\\(1) $\mu=\pm 3, r=s=3$ and $G \cong K_{3,3}$;
	\\(2) $1 \neq \mu \in \mathbb{N}_{+}, r=s$ and $G$ is an $r$-regular graph of order $\left(\mu^{4}+10 \mu^{3}+27 \mu^{2}+10 \mu\right) / 4$, where $r=\mu(\mu+1)(\mu+4) / 2$.
	\\(3) $\mu=1$, $s=5$ and either $G\cong Sch_{10}$ or $G$ is isomorphic to one of the eleven induced regular subgraphs of $Sch_{10}$. 
	\\(4) $\mu=-1$, $r\equiv -1\ ({\rm mod} \ 2s-1)$ and $G\cong G'(r)$ (see \cite{yuan} for specific definitions).
\end{theo}

In this section, we consider the general case. We prove that there is no regular graph $G$ with $K_{t,s}\ (s\geq t\geq 1)$ as a star complement for $\mu=-t$, characterize the graph $G$ when $\mu=r$, $\mu=-1$, and the case with all vertices in $X$ of type $(0,b)$ for $\mu\notin \{-t,r,-1\}$. Furthermore, we propose a question for further research.

\begin{prop}\label{u=-t}
	There is not an $r$-regular graph $G$ with $K_{t,s}\ (s\geq t\geq 1)$ as a star complement for $\mu=-t$.
\end{prop}
\begin{Proof}
	Let $\mu = -t$. Since $\mu^2 \neq ts$, we have $s\neq t$ and then $s>t$. Let $u\in X$ be a vertex of type $(a,b)$, thus $(a,b)\neq (0,0)$ and $0\leq a\leq t$, $0\leq b\leq s$. 
	
	If $t=1$, from Theorem 2.2 of \cite{9}, there is no $r$-regular graph $G$ with $K_{1,s}$ as a star complement for $\mu=-1$.
	
	If $t\geq 2$, by Lemma \ref{lem25}, we know that $\mu=-t$ is a non-main eigenvalue of $G$, and by (\ref{eq25}), we have
	\begin{equation}\label{eq32'}
	t(t-s)(a-t)=0.
	\end{equation}
	
	Since $s> t$ and $t\geq 2$, (\ref{eq32'}) implies that $a=t$, and thus $s(b-t^2)=b^2-tb-t^3+t^2$ by (\ref{eq26}). 
	
	If $b=t^2$, then $t^4-2t^3+t^2=0$, thus $t=0$ or $1$, a contradiction. 
	
	If $b\neq t^2$, then $s=\frac{b^2-tb-t^3+t^2}{b-t^2}=b-t^2+\frac{t^2(t-1)^2}{b-t^2}+2t^2-t$. If $b<t^2$, then $s\leq -2\sqrt{t^2(t-1)^2}+2t^2-t=t$, which contradicts with $s>t$.
	Thus $b> t^2$ and  $s-(b+t-1)=\frac{b(t-1)^2}{b-t^2} >0$. Considering degrees, we have
	$$d_G(v_1)=d_G(v_2)=\cdots=d_G(v_t)=s+\left|X\right|,$$
	and
	$$d_G(u)\leq a+b+\left|X\right|-1=b+t-1+|X|,\ u\in X.$$
	Hence, $d_G(v_1)=d_G(v_2)=\cdots=d_G(v_t)>d_G(u)$ which contradicts to the regularity of $G$. 
	
	Combining the above arguments, there is not an $r$-regular graph $G$ with $K_{t,s}\ (s\geq t\geq 1)$ as a star complement for $\mu=-t$.$\quad \square$
\end{Proof}

\begin{theo}\label{u=r}
	If the $r$-regular graph $G$ has $K_{t,s}\ (s\geq t\geq 1)$ as a star complement for an eigenvalue $\mu=r$, then $s=t=1$, $G\cong C_3$ or $r=s=t+1$, $G\cong K_{t+1,t+1}$.
\end{theo}
\begin{Proof}
	Since $\mu$ is not an eigenvalue of $H\cong K_{t,s}\ (s\geq t\geq 1)$, we have $\mu \neq 0$ and $\mu^2 \neq ts$. By Lemma \ref{lem23}, $V(K_{t,s})$ is a location-dominating set, and so $G$ is connected.
	
	By $G$ is $r$-regular and connected, $\mu=r$, we know $k=1$ and then $\left|X\right|=1$. Since $G$ is regular, we have $d_G(v_1)=d_G(v_2)=\cdots=d_G(v_t)$. Let $X=\left\{u\right\}$. Then either $u\sim v_1$, $u\sim v_2$, \dots,$u\sim v_t$, or $u\nsim v_1$, $u\nsim v_2$,\dots, $u\nsim v_t$. 
	
	If $u\sim v_1$, $u\sim v_2$, \dots,$u\sim v_t$, then $d_G(v_1)=d_G(v_2)=\cdots=d_G(v_t)=s+1,$ which implies that $d_G(u)=s+1$. It follows that the vertex $u$ is adjacent to $s-t+1$($\geq 1$) vertices of $W$, and thus $t=1$ by $d_G(w_1)=d_G(w_2)=\cdots=d_G(w_s)$. Since $d_G(w_1)=t+1=d_G(v_1)$, we have $s=t=1$ and $G\cong C_3$.
	
	If $u\nsim v_1$, $u\nsim v_2$, \dots, $u\nsim v_t$, in view of the regularity, we have $d_G(u)=d_G(v_1)=d_G(v_2)=\cdots=d_G(v_t)=s,$ and then $d_G(w_1)=d_G(w_2)=\cdots=d_G(w_s)=t+1$. Hence we have $s=r=t+1$ and $G\cong K_{t+1,t+1}$.$\square$
\end{Proof}

Let $H\cong K_{t,s}$, $(V, W)$ be a bipartition of the graph $K_{t,s}$ with $V=\left\{ v_1,v_2,\dots, v_t \right\}$, $W=\left\{ w_1,w_2,\dots,w_s \right\}$. We obtain an $r$-regular graph $G(r)$ with $V(G(r))=X\cup V(H)$, $X=V_1 \cup \cdots \cup V_t \cup W_1 \cup \cdots \cup W_s$ where $V_i$ is the set of vertices of type $(1,s)$ adjacent to $v_i\in V$ with $|V_i|=(r+1)(s-1)/(ts-1)-1$, $V_i$ induces a clique for $1\leq i \leq t$, $W_i$ is the set of vertices of type $(t,1)$ adjacent to $w_i\in W$ with $|W_i|=(r+1)(t-1)/(ts-1)-1$, $W_i$ induces a clique for $1\leq i \leq s$. For any $i$, $j$, each vertex in $V_i$ is adjacent to all vertices in $W_j$. 




The greatest common divisor of $a$ and $b$ is denoted by $\gcd(a,b)$. For $\mu=-1$, we have the following theorem.

\begin{theo}\label{u=-1}
	If $G$ is an $r$-regular graph with $H=K_{t,s}\ (s\geq t\geq 2)$ as a star complement for an eigenvalue $-1$, then $r\equiv -1\ ( {\rm mod} \frac{ts-1}{\gcd(s-1,t-1)})$ and $G\cong G(r)$.
\end{theo}
\begin{Proof}
	Since $K_{t,s}$ is connected and $V(K_{t,s})$ is a dominating set (see Lemma \ref{lem23}), we know $G$ is connected. Let $H\cong K_{t,s}$, $(V, W)$ be a bipartition of the graph $K_{t,s}$ as above. Let $u\in X$ be a vertex of type $(a,b)$, thus $(a,b)\neq (0,0)$ and $0\leq a\leq t$, $0\leq b \leq s$. Let $\mu=-1$ in (\ref{eq26}), so that
	\begin{equation}\label{1u1}
	(1-ts)(a+b-1)+a^2s+tb^2-2ab=0.
	\end{equation}

	By Lemma \ref{lem25}, we know that $\mu=-1$ is a non-main eigenvalue of $G$, thus from (\ref{eq25}), we have
	\begin{equation}\label{1u2}
	1-ts+a(s-1)+b(t-1)=0.
	\end{equation}

	Combining (\ref{1u1}) and (\ref{1u2}), if $a=1$, then $b=s$; if $a=t$, then $b=1$; if $1<a<t$, then $b=a-t+1,\ s=2-t$ or $b=\frac{a}{t},\  s=\frac{1}{t}$. It is obvious that $s=2-t$ or $s=\frac{1}{t}$ contradicts with $s\in \mathbb{Z}_+$. Therefore, the possible types of vertices in $X$ are $(1,s)$, $(t,1)$, and the feasible solution of (\ref{eq27}) are shown in Table \ref{solution}.
	\begin{table}[h]
		\centering
		\begin{tabular}{llll}
			\hline
			$(a,b)$    & $(c,d)$  & $a_{uv}$ & $\rho_{uv}$ \\
			\hline
			$(1,s)$   &  $(1,s)$ & $0$ & $s$ \\
			$(1,s)$  &  $(1,s)$ & $1$ & $s+1$ \\
			$(t,1)$ & $(t,1)$ & $0$ & $t$\\
			$(t,1)$ & $(t,1)$& $1$ & $t+1$ \\
			$(1,s)$ & $(t,1)$ & $1$ & $2$ \\
			\hline
		\end{tabular}
		\caption{The feasible solution of (\ref{eq27})}
		\label{solution}
	\end{table}
	
	We observe that when $u,v$ are of different types, they must be adjacent; when $u,v$ are of the same type, $u\sim v$ if and only if they have the same $H$-neighbourhoods, and thus $u,v$ are co-duplicate vertices. We can add arbitrarily many co-duplicate vertices when constructing graphs with a prescribed star complement for $-1$.
	
	Now we partition the vertices in $X$. Let $V_i$ be the set of vertices of type $(1,s)$ in $X$ adjacent to $v_i\in V$, $W_i$ be the set of vertices of type $(t,1)$ in $X$ adjacent to $w_i\in W$. Then any two vertices in $V_i\ (W_i)$ are co-duplicate vertices. We do not exclude the possibility that some of the sets $V_i$, $W_i$ are empty. Then for any $v_i\in V$, we have $d_G(v_i)=s+|V_i|+\sum\limits_{i=1}^s |W_i|;$ and for any $w_i\in W$, we have $d_G(w_i)=t+\sum\limits_{i=1}^t |V_i|+|W_i|.$
	
	Since $G$ is $r$-regular, we have $|V_1|=|V_2|=\cdots=|V_t|$ by $d_G(v_1)=d_G(v_2)=\cdots =d_G(v_t)$ and $|W_1|=|W_2|=\cdots=|W_s|$ by $d_G(w_1)=d_G(w_2)=\cdots =d_G(w_s)$. Then we have 
	$$r=d_G(v_1)=s+|V_1|+s\cdot |W_1| \ \mbox{and} \
	r=d_G(w_1)=t+t\cdot |V_1|+|W_1|.$$
	It turns out that
	$$|V_1|=\frac{(s-1)(r+1)}{ts-1}-1,\ |W_1|=\frac{(t-1)(r+1)}{ts-1}-1.$$
	Since $|V_1|\in \mathbb{N}$, $|W_1|\in \mathbb{N}$ and 
	$$\gcd(t-1,ts-1)=\gcd(s-1,ts-1)=\gcd(t-1,s-1),$$
	we have $r\equiv -1\ ({\rm mod} \ \frac{ts-1}{\gcd(s-1,t-1)})$. Consequently we obtain an $r$-regular graph $G(r)$.$\quad \square$
\end{Proof}

\begin{remark}
	Note that if $\ V_i^{*}=V_i\cup \left\{v_i\right\},\ v_i\in V$ and $W_i^*=W_i\cup \left\{w_i\right\}, \ w_i\in W$, then each of sets $V_i^*$, $W_i^*$ induces a clique in $G(r)$.
\end{remark}

Next, we consider the case $\mu \notin \left\{-1,-t,r\right\}$. The following lemma lists all possible types of vertices in $X$.
\begin{lem}\label{type}
	Let $G$ be a graph with $H=K_{t,s}\ (s\geq t\geq 1)$ as a star complement for $\mu$. If $\mu$ is a non-main eigenvalue of $G$ and $\mu \notin \left\{-1,-t\right\}$, then the possible types of vertices in the star set $X$ are $(a, \frac{\mu^3+t\mu^2-ta+a^2}{\mu +a})$, where $0\leq a \leq t-1$ and $a\neq -\mu $.
\end{lem}
\begin{Proof}
  Let $u \in X$ be a vertex of type $(a,b)$, thus $(a,b)\neq (0,0)$ and $0 \leq a \leq t$, $0\leq b\leq s$. By (\ref{eq25}), (\ref{eq26}) and $\mu \notin\{0,-1,-t\}$, we have: (1) $a=t$, $b=-\mu$, $s=\frac{\mu^2}{t}$; (2) $a=-\mu$, $b=\frac{\mu^2}{t}$, $s=\frac{\mu^2}{t}$; (3) $0\leq a \leq t-1$, $a\neq -\mu$, $\left\{\begin{matrix*}[l] b=\frac{-a\mu}{t},\\ 
  s=\frac{\mu^2}{t},	
  \end{matrix*}\right.$  or
   $\left\{\begin{matrix*}[l]
  	b=\frac{\mu^3+t\mu^2-ta+a^2}{\mu +a},\\
  	s=\frac{\mu^4+(2t+1)\mu^3+(2a+t^2)\mu^2+(2a^2-at)\mu+a^2t-at^2}{(t-a)\mu+at-a^2}.
  	\end{matrix*}\right.$
  	Since $\mu$ is not an eigenvalue of $H$, we have $\mu^2 \neq ts$. Thus the possible types of vertices in the star set $X$ are $(a, \frac{\mu^3+t\mu^2-ta+a^2}{\mu +a})$.$\qquad\square$
\end{Proof}

For $\mu\notin \{-1,-t,r\}$, we consider the case $a=0$ in the following by Lemma \ref{type}.

\begin{theo}\label{k0b}
	If the $r$-regular graph $G$ has $K_{t,s}\ (s\geq t\geq 1)$ as a star complement for $\mu \notin \left\{-1,-t,r\right\}$ and all vertices in $X$ are of type $(0,b)$, then one of the following holds:
	\\{\rm(1)} $\mu=-r, r=s=t+1$ and $G \cong K_{t+1,t+1}$;
	\\{\rm(2)} $\mu=\frac{1}{2}(-1\pm \sqrt{5}), t=1, r=s=2$ and $G$ is a $5$-cycle;
	\\{\rm(3)} $\mu \in \mathbb{N}_{+}, r=s$ and $G$ is an $r$-regular graph of order $\mu \left( \mu+2t+1\right) \left(\mu^{2}+2t \mu+\mu+t^2-t \right) / t^2$, where $r=\mu \left(\mu^{2}+2t \mu+\mu+t^2 \right) / t$.
\end{theo}
\begin{Proof}
	Since $\mu$ is not an eigenvalue of $K_{t,s}$, we have $\mu \neq 0$ and $\mu^2 \neq ts$. By Lemma \ref{lem23}, $V(K_{t,s})$ is a location-dominating set, and so $G$ is connected. Then by $a=0$ and Lemma \ref{type}, we have
	\begin{equation*}
	\left\{\begin{matrix}
	b=\mu^2+t\mu,\\
	s=\mu \left(\mu^{2}+2t \mu+\mu+t^2 \right) / t.
	\end{matrix}\right.
	\end{equation*}
	
	Now we consider $H=K_{t,\mu \left(\mu^{2}+2t \mu+\mu+t^2 \right) / t}$ and all vertices in $X$ are of type $(0, \mu^2+t\mu)$. Then $r=s$. Counting the edges between $X$ and $V(H)$, we have
	$\left|X\right|(\mu^2+t\mu)=s(r-t).$
	Thus
	\begin{equation}\label{eq43}
	\left|X\right|=\frac{s(r-t)}{(\mu^2+t\mu)}=\frac{1}{t^2}(\mu^2+2t\mu+\mu+t^2)(\mu^2+t\mu+\mu-t).
	\end{equation}
	
\noindent	\textbf{Case 1}: $\left|X\right|=1$.
	
	 Then $(\mu+t+1)\left(\mu^3+(2t+1)\mu^2+(t^2-t)\mu-t^2 \right)=0$ from (\ref{eq43}). When $\mu^3+(2t+1)\mu^2+(t^2-t)\mu-t^2=0$, then $\mu \notin \mathbb{Z}$ and thus $r=\frac{\mu}{t} \left(\mu^{2}+2t \mu+\mu+t^2 \right)-\frac{1}{t}(\mu^3+(2t+1)\mu^2+(t^2-t)\mu-t^2)=\mu+t\notin \mathbb{Z}$, a contradiction. When $\mu=-(t+1)$, then $r=s=t+1$ and $G\cong K_{t+1,t+1}$.
	
\noindent	\textbf{Case 2}: $\left|X\right|\geq 2$. 
	
	We apply the compatibility condition (\ref{eq27}) to vertices $u,v$ in $X$, we find that
	\begin{equation}\label{trho}
	\rho_{uv}=\left\{ \begin{matrix}
	t\mu, & u\nsim v,\\
	(t-1)\mu, & u\sim v.
	\end{matrix}\right.	
	\end{equation}
	
	If $t=1$ and $X$ induces a clique then $|X|-1=r-\mu^2-\mu$, whence $(\mu+1)(\mu+2)(\mu^2+\mu-1)=0$. Thus either $\mu=-2$, $r=s=t+1=2$ and $G\cong K_{2,2}$ which belongs to case (1), or $\mu=\frac{1}{2}(-1\pm \sqrt{5})$ and we have case (2) (\cite{9}).
	
	Otherwise, it follows from (\ref{trho}) that $\mu \in \mathbb{N}_+$ and $G$ is $\mu \left(\mu^{2}+2t \mu+\mu+t^2 \right) / t$-regular of order $\mu( \mu+2t+1) (\mu^{2}+2t \mu+\mu+t^2-t )/t^2$ with $K_{t,\mu \left(\mu^{2}+2t \mu+\mu+t^2 \right) / t}$ as a star complement for $\mu$.$\quad \square$
\end{Proof}

In \cite{9}, Rowlinson gave a lemma to determine whether a connected $r$-regular graph with $K_{1,s}$ as a star complement is a strongly regular graph. Now we extend it to $K_{t,s}$.
\begin{lem}\label{lem31}
	Let $G$ be a connected $r$-regular graph with $\mu(\neq r)$ as an eigenvalue of multiplicity $k$. Suppose that $\left |V(G)\right|=k+t+s$. If $k+t+s-1>r$, then
	$$(k+t+s)r-r^2-k\mu^2-\frac{(k\mu +r)^2}{s+t-1}\geq 0,$$
	with equality if and only if $G$ is strongly regular.
\end{lem}
\begin{Proof}
	Since $k+t+s-1>r$, neither $G$ nor $\overline{G}$ is complete. Let $\theta_1,\dots,\theta _{s+t-1}$ be the eigenvalue of $G$ other than $\mu$ and $r$. We have
	$$ \sum _{i=1}^{s+t-1} \theta_i+k\mu+r=0 \mbox{ and } \sum _{i=1}^{s+t-1} \theta_{i}^2+k\mu^2+r^2=(k+t+s)r.$$
	It follows that if $\overline{\theta}=\frac{1}{s+t-1} \sum\limits _{i=1}^{s+t-1}\theta_i$, then
	$$\sum\limits _{i=1}^{s+t-1}(\theta_i-\overline{\theta} )^2=\sum\limits _{i=1}^{s+t-1} \theta_i^2-(s+t-1)\overline{\theta}^2=(k+t+s)r-r^2-k\mu^2-\frac{(k\mu +r)^2}{s+t-1}\geq 0.$$
	Equality holds if and only if $\theta_i=\overline{\theta}\ (i\in [s+t-1])$, equivalently $G$ has just three distinct eigenvalue. By \cite[Theorem 1.2.20]{4}, a non-complete connected regular graph is strongly regular if and only if it has exactly three distinct eigenvalues. The proof is completed.$ \quad \square$
\end{Proof}
\begin{remark}
	From Lemma \ref{lem31}, we know that the $r$-regular graph $G$ in {\rm (3)} of Theorem \ref{k0b} is strongly regular when $t=1$, and the graph $G$ isn't strongly regular when $t\geq 2$.
\end{remark}

For the cases $1\leq a \leq t-1$, we cannot give a characterization. Thus we propose a question for further research.

\begin{que}\label{que}
	Let $s\geq t\geq 1$, $\mu\notin \{-1,-t,r\}$ and $1\leq a\leq t-1$. Can we give a characterization of the $r$-regular graphs with $K_{t,s}$ as a star complement?
\end{que}

\section{Regular graphs with $K_{3,s}$ as a star complement}\label{regular}
\hspace{1.5em}
In this section, we completely solve Question \ref{que} when $t=3$. Since the cases when $\mu \in \left\{-1,-3,r\right\}$ have been solved in Section \ref{kts}, we consider the cases $\mu \notin \left\{-1,-3,r\right\}$ in the following. 

Let $G$ be an $r$-regular graph with $H=K_{3,s}\ (s\geq 3)$ as a star complement for $\mu \notin \left\{-1,-3,r\right\}$. Then $\mu \neq 0$ and $\mu^2\neq 3s$ for $\mu$ is not an eigenvalue of $K_{3,s}$.
By Lemma \ref{type}, the possible types of vertices in $X$ are shown in Table \ref{posstype}:
\begin{table}[h]
\centering
\begin{tabularx}{10cm}{lll}
\hline
Type & $(a,b)$ & s \\
\hline  \\[-4mm]
I  & $(0,\mu^2+3\mu)$ & $\mu (\mu^2+7\mu+9)/3$ \\ [2pt]
II  & $(1,\mu^2+2\mu-2)$ & $(\mu+2)(\mu^2+4\mu-3)/2$ \\ [2pt]
III & \large  $(2,\frac{\mu ^3+3\mu^2-2}{\mu+2})$  &  \large $\frac{\mu^4+7\mu^3+13\mu^2+2\mu-6}{\mu+2}$ \\[2pt]
\hline
\end{tabularx}
\caption{The possible types of vertices in $X$}
\label{posstype}
\end{table}

\begin{lem}\label{lem34}
Let $G$ be a graph with $H=K_{3,s}\ (s\geq 3)$ as a star complement for $\mu$. If $\mu$ is a non-main eigenvalue of $G$ and $\mu \notin \left\{-1,-3\right\}$, then all the vertices in the star set $X$ are of the same type, say, Type I, Type II or Type III.
\end{lem}
\begin{Proof}
Now we show it is impossible that there are two or three types of vertices in $X$.

\noindent \textbf{Case 1.} There exist vertices of Type I and Type III in $X$.

From Table \ref{posstype}, we have
$$s=\frac{1}{3}\mu (\mu^2+7\mu+9)=\frac{\mu^4+7\mu^3+13\mu^2+2\mu-6}{\mu+2}$$
with solution $\mu=-3$, $-1$ or $1$. Since $\mu \notin \{-1,-3\}$, we have $\mu=1$, and thus $s=17/3$, a contraction.

\noindent \textbf{Case 2.} There exist vertices of Type II and Type III in $X$.

From Table \ref{posstype}, we have
$$s=\frac{1}{2}(\mu+2)(\mu^2+4\mu-3)=\frac{\mu^4+7\mu^3+13\mu^2+2\mu-6}{\mu+2}$$
with solution $\mu=-3$ or $0$, a contraction.

\noindent \textbf{Case 3.} There exist vertices of Type I and Type II in $X$.

By Case 1 and Case 2, we know there does not exist vertices of Type III in $X$.

By Table \ref{posstype}, we have
$$s=\frac{1}{3}\mu (\mu^2+7\mu+9)=\frac{1}{2}(\mu+2)(\mu^2+4\mu-3)$$
with solution $\mu=-3$ or $2$. Since $\mu \neq -3$, we have $\mu=2$, and thus $s=18$. The vertices in $X$ are of type $(0,10)$, $(1,6)$. From (\ref{eq27}), Table \ref{table2} is obtained.
\begin{table}[h]
    \centering
    \begin{tabular}{llcc}
    \hline
      $(a,b)$   & $(c,d)$ & $a_{uv}$ & $\rho_{uv}$\\
      \hline
      $(0,10)$   & $(1,6)$ & $0$ & $4$\\
      $(0,10)$ & $(1,6)$ & $1$ & $2$\\
       $(0,10)$ & $(0,10)$ & $0$ & $2$\\
       $(0,10)$  & $(0,10)$ & $1$ & $0$\\
        $(1,6)$  & $(1,6)$ & $0$ & $-11/5$\\
          $(1,6)$ & $(1,6)$ & $1$ & $-21/5$\\
          \hline
    \end{tabular}
\caption{The possible adjacency of vertices in $X$}
\label{table2}
\end{table}

For any two vertices of type $(1,6)$, $(1,6)$ in $X$, $\rho_{uv} \notin \mathbb{Z}$, so there is at most one vertex of type $(1,6)$ in $X$. For any two vertices of type $(0,10)$, $(0,10)$ in $X$, $\rho_{uv} =0$ or $2$. Since $s=18$, there are at most two vertices of type $(0,10)$ in $X$, and $\rho_{uv}=2$ if there are two vertices of type $(0,10)$. Therefore, $X$ has at most three vertices.

Let $V(K_{3,18})=V\cup W$ with $\left|V\right|=3$, $\left|W\right|=18$. For any $v_i\in V$ and $w_i\in W$, we have $d_G(v_i)\geq 18$ and $d_G(w_i)\leq 3+3=6$, which contradicts with the regularity of $G$.

The proof is completed.$\quad \square$
\end{Proof}

\begin{figure}[h]
	\centering
	\subfloat[$G_1$]{
		\begin{minipage}{180pt}
			\centering
			\includegraphics[scale=0.12]{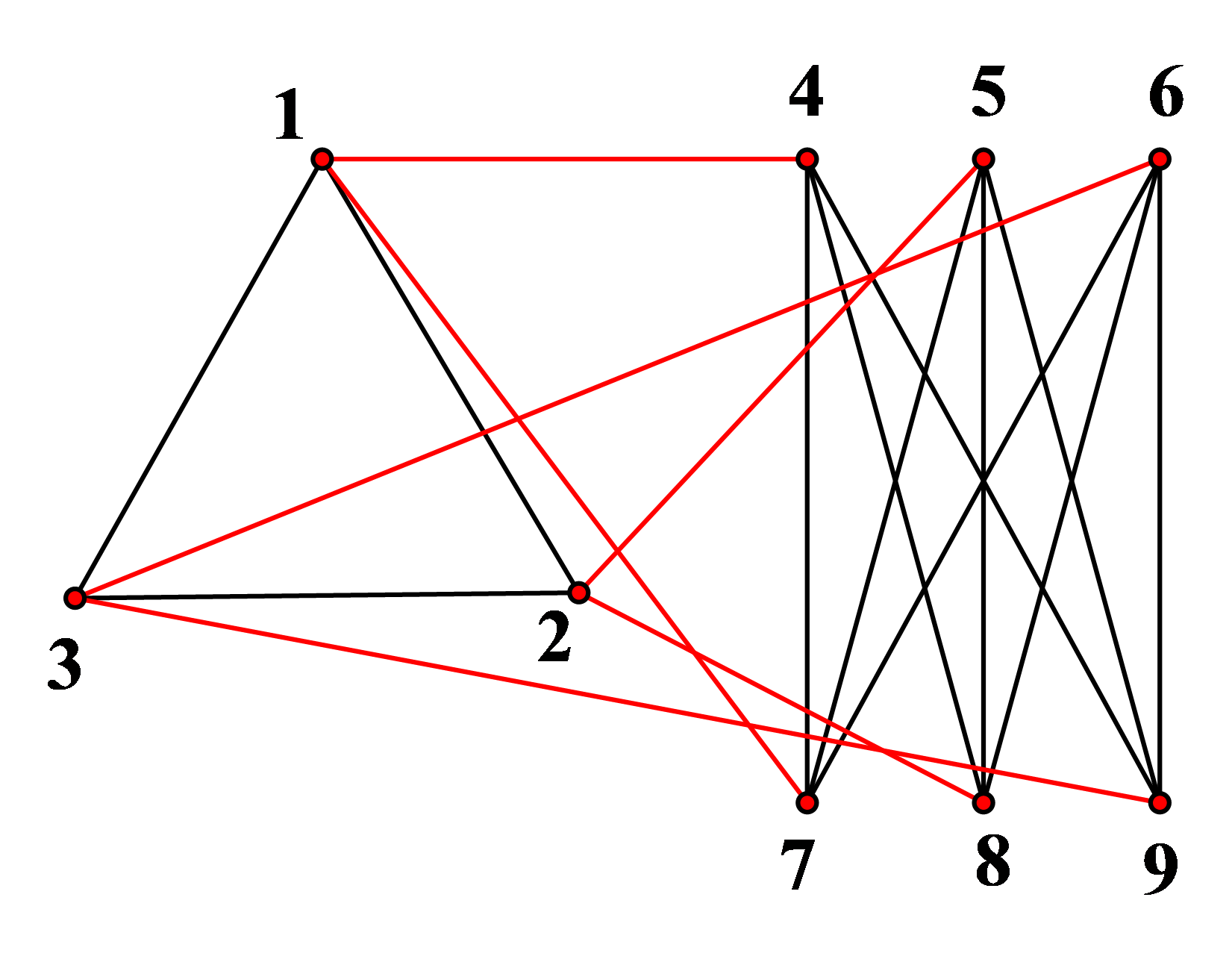}
	\end{minipage}}
	\qquad
	\subfloat[$G_2$]{
		\begin{minipage}{260pt}
			\centering
			\includegraphics[scale=0.1]{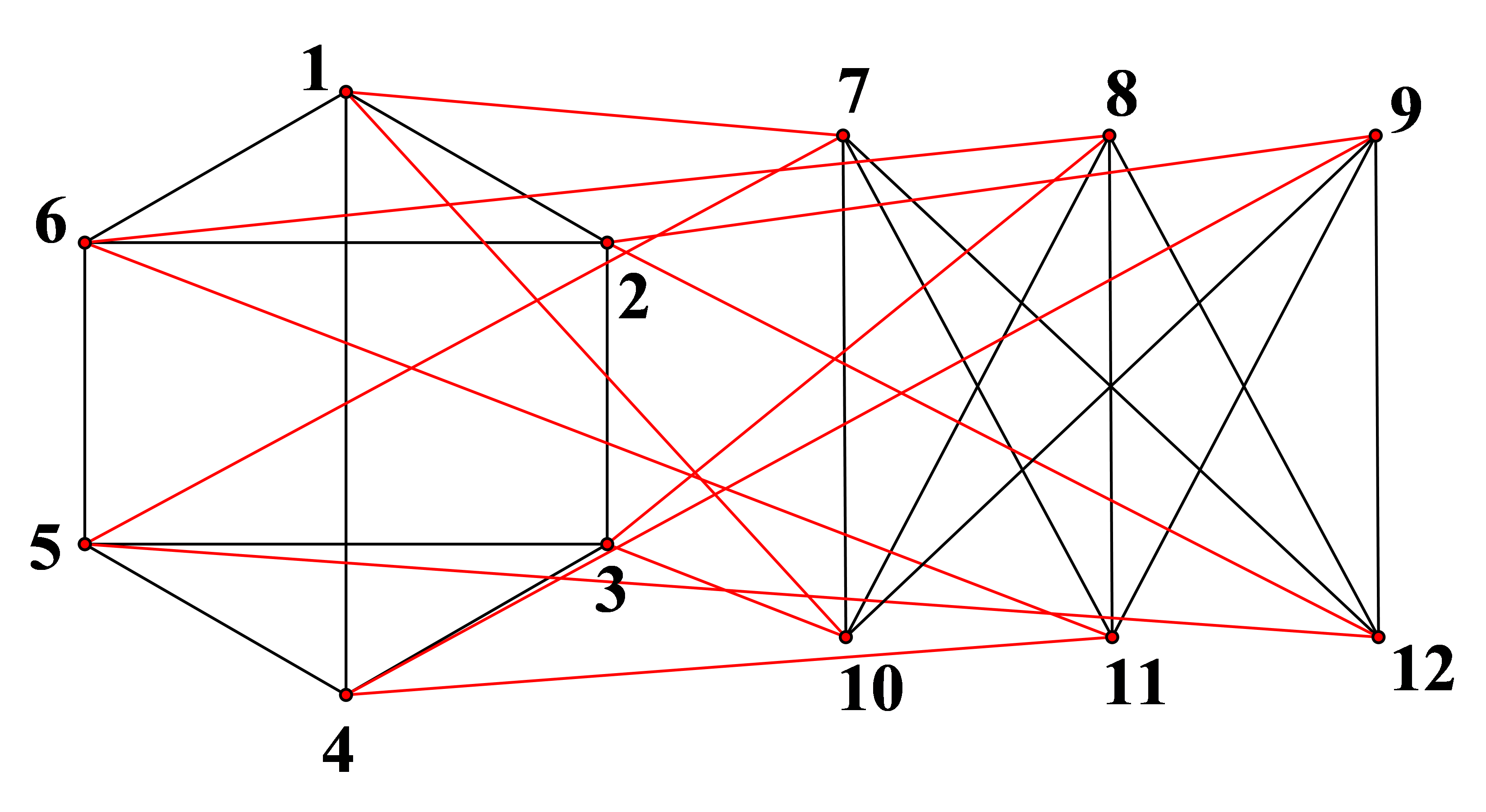}
	\end{minipage}}
	\qquad
	\subfloat[$G_3$]{
		\begin{minipage}{450pt}
			\centering
			\includegraphics[scale=0.1]{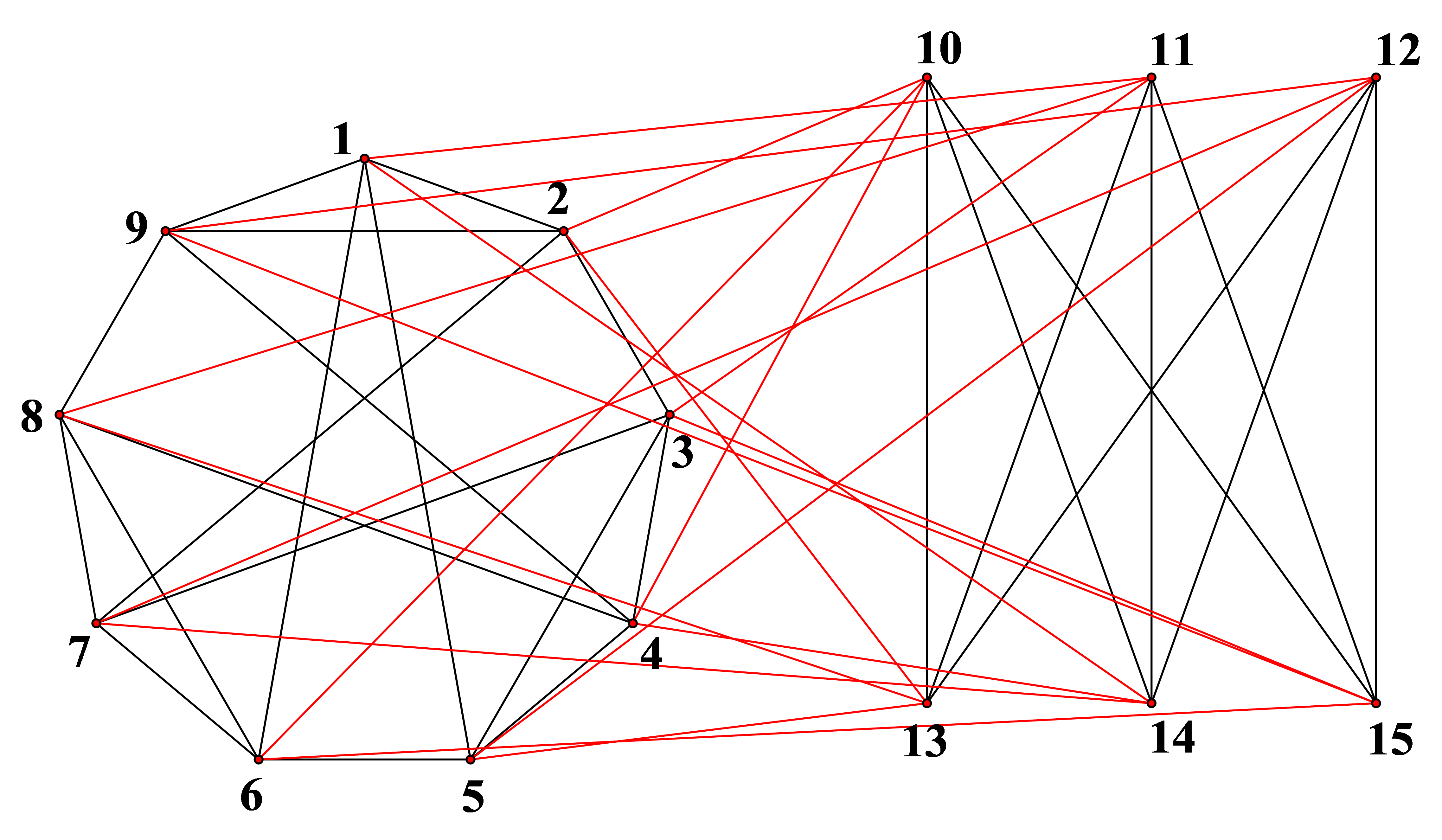}
	\end{minipage}}
	\caption{Regular graphs $G_1$, $G_2$, $G_3$ of Theorem \ref{thm35}}
	\label{mugraph}
\end{figure}

Now we define three special graphs $G_1, G_2$ and $G_3$ (see Figure \ref{mugraph}) as follows. Clearly, the graph $G_1$ is a 4-regular graph of order 9, its spectrum is [$-3,-2^2,0^2,1^3,4$];
$G_2$ is a 5-regular graph of order 12, its spectrum is [$-3^3,-1^2,1^6,5$];
$G_3$ is a 6-regular graph of order 15, its spectrum is [$-3^5,1^9,6$]. By Lemma \ref{lem31}, we know that $G_1$ and $G_2$ isn't strongly regular while $G_3$ is strongly regular with parameters $(15,6,1,3)$.

\begin{theo}\label{thm35}
Let $s\geq 3$, $G$ be an $r$-regular graph with $K_{3,s}$ as a star complement for an eigenvalue $\mu$. Then one of the following holds:\\
{\rm(1)} $\mu=-1$, $r\equiv -1\ ({\rm mod }\ \frac{3s-1}{(s-1,2)})$ and $G\cong G(r)$, where $G(r)$ is defined in Theorem \ref{u=-1};\\
{\rm(2)} $\mu=\pm 4$, $r=s=4$ and $G\cong K_{4,4}$;\\
{\rm(3)} $\mu \in \mathbb{N}_+$, $r=s$ and $G$ is an $r$-regular graph of order $\mu(\mu+7)(\mu+6)(\mu+1)/9$, where $r=\mu(\mu^2+7\mu+9)/3$;\\
{\rm(4)} $\mu=1$, $s=3$, and $G\cong G_1\ (r=4)$, $G\cong G_2\ (r=5)$ or $G\cong G_3\ (r=6)$ {\rm(see Figure \ref{mugraph})}.
\end{theo}
\begin{Proof}
By Theorem \ref{u=-1} and $\mu=-1$, we have (1) holds.

By Theorem \ref{u=r} and $\mu=r$, we have $s=r=4$ and $G\cong K_{4,4}$.

If $\mu\notin \{-1,r\}$, then $\mu$ is non-main. From Proposition \ref{u=-t}, we know $\mu \neq -3$. Since $\mu$ is not an eigenvalue of $H\cong K_{3,s}$, we have $\mu \neq 0$ and $\mu^2 \neq 3s$. By Lemma \ref{lem23}, $V(K_{3,s})$ is a location-dominating set, and so $G$ is connected. Let $V(K_{3,s})=V\cup W$ with $\left|V\right|=3$, $\left|W\right|=s$. Denote the vertices in $V$ and $W$ by $v_1,v_2,v_3$ and $w_1,w_2,\dots,w_s$, respectively. In the following, we suppose that $\mu \neq \left\{-1,r,-3,0\right\}$ and $\mu^2\neq 3s$. From Table \ref{posstype}, there are only three possible types of vertex in $X$. By Lemma \ref{lem34}, we know all the vertices in $X$ are of the same type, and we consider the following three cases:

\noindent\textbf{Case 1.} The vertices in $X$ are of Type I. 

Then (1) or (3) of Theorem \ref{k0b} holds by $s\geq 3$, say, either $\mu=-4$, $r=s=4$ and $G\cong K_{4,4}$, or $\mu \in \mathbb{N}_+$ and $G$ is $\mu(\mu^2+7\mu+9)/3$-regular of order $\mu(\mu+7)(\mu+6)(\mu+1)/9$ with $K_{3,(\mu^3+7\mu^2+9\mu)/3}$ as a star complement for $\mu$. 

Combining the above case of $\mu=r$, result (2) or (3) holds.

\noindent\textbf{Case 2.} The vertices in $X$ are of Type II. 

Then $H=K_{3,(\mu+2)(\mu^2+4\mu-3)/2}$ and all vertices in $X$ are of type $(1,\mu^2+2\mu-2)$. Applying the compatibility condition (\ref{eq27}) to vertices $u,v$ of $X$, since $\mu\neq -3$, we find that
\begin{equation}\label{rho2}
    \rho_{uv}=\left\{
    \begin{matrix}
    \mu-1, & u\sim v,\\
    2\mu-1, & u\nsim v.
    \end{matrix}
    \right.
\end{equation}

By regularity of $G$, we have $d_G(v_1)=d_G(v_2)=d_G(v_3)$. This implies the vertices in $X$ are equally divided into three parts, and each vertex in $V$ is adjacent to every vertex of one part, so
\begin{equation}\label{rx}
   r=s+\left|X\right|/3.
\end{equation}

Now we compute the edges between $X$ and $V(H)$ by two ways, and we have
\begin{equation}\label{rx2}
\left|X\right|(\mu^2+2\mu-1)=3(r-s)+s(r-3).
\end{equation}
By (\ref{rx}), (\ref{rx2}) and $s=(\mu+2)(\mu^2+4\mu-3)/2$, we have $(\mu+3)(\mu-1)(\mu-2) \left|X\right|= -\frac{3}{2}(\mu+2)(\mu^2+4\mu-3)(\mu+3)(\mu-1)(\mu+4).$
Since $\mu \neq -3$, we consider the following three subcases. 

\noindent\textbf{Subcase 2.1.} $\mu\neq 1,2$.

Then
$$\left|X\right|=-\frac{3}{2}(\mu+2)(\mu^2+4\mu-3)\frac{\mu+4}{\mu-2}=-3s(1+\frac{6}{\mu-2}).$$
Since $\left|X\right|\in \mathbb{Z}$ and $s\in \mathbb{Z}$, we have $\mu \in \mathbb{Q}$. Notice that $\mu$ is an algebraic integer, then $\mu \in \mathbb{Z}$. From (\ref{rho2}), we know that $\mu\in\mathbb{N}_+\setminus \left\{1,2\right\}$. Thus $\left|X\right|<0$, a contradiction.

\noindent\textbf{Subcase 2.2.} $\mu=2$.

 Then $s=\frac{(\mu+2)(\mu^2+4\mu-3)}{2}=18$, $H=K_{3,18}$ and all vertices in $X$ are of type $(1,6)$.
By (\ref{rx}) and (\ref{rx2}), we obtain $0=270$, it is a contradiction.

\noindent\textbf{Subcase 2.3.} $\mu=1$.

Then $s=3$, $H=K_{3,3}$ with $V=\left\{ v_1,v_2, v_3 \right\}$, $W=\left\{ w_1,w_2,w_3 \right\}$ and all vertices in $X$ are of type $(1,1)$. Thus $\left|X\right|\leq \binom{3}{1}\binom{3}{1}=9$. From (\ref{rho2}), we have
\begin{equation}\label{u1rho}
    \rho_{uv}=\left\{
    \begin{matrix}
    0, & u\sim v,\\
    1, & u\nsim v.
    \end{matrix}
    \right.
\end{equation}
Since $G$ and $H$ are regular, $X$ induces a $r'$-regular graph, denoted by $G[X]$. 
\vskip 0.2cm

\noindent \textbf{Claim 1.} 
	The graph $G[X]$ cannot contain $S_1$, $S_2$ or $S_3$ as an induced subgraph (see Figure \ref{subgra}). 
		
\begin{Proof}
	If $G[X]$ contains $S_1$ as an induced subgraph,  without loss of generality, we suppose that $u_1\sim v_1, u_1\sim w_1$. Then by (\ref{u1rho}), $\rho_{u_1u_2}=\rho_{u_1u_3}=\rho_{u_1u_4}=0$, and thus vertices $u_2$, $u_3$ and $u_4$ are adjacent to one vertex in $V\setminus \left\{v_1\right\}$ and one vertex in $W\setminus \left\{w_1\right\}$ such that $\rho_{u_2u_3}=\rho_{u_2u_4}=\rho_{u_3u_4}=1$, it is impossible.
	
	If $G[X]$ contains $S_2$ as an induced subgraph, since the vertices $u_5$, $u_6$ and $u_7$ are adjacent in pairs, by (\ref{u1rho}), without loss of generality, we can suppose that $u_5\sim v_1, u_5\sim w_1, u_6\sim v_2, u_6\sim w_2, u_7\sim v_3, u_7\sim w_3$. Since $u_8\sim u_6, u_8\sim u_7$, by (\ref{u1rho}), vertex $u_8$ is adjacent to vertices in $V(K_{3,3})\setminus \left\{v_2,w_2,v_3,w_3\right\}$. Thus the $H$-neighbourhood of $u_5$ and $u_8$ is the same, a contradiction.
	
	Similar to the proof of $S_2$, it's obvious that $G[X]$ cannot contain $S_3$ as an induced subgraph.$\qquad \square$
\end{Proof}
\vskip 0.2cm

By $G[X]$ is $r'$-regular and (\ref{rx}), for any $u\in X$, $d_G(u)=2+r'=r=3+\left|X\right|/3$. Then $2\leq r'= 1+\left|X\right|/3 \leq 4$ by $|X|\leq 9$.

If $r'=2$, then $\left|X\right|=3$, $G[X]\cong C_3$, $G\cong G_1$ by (\ref{u1rho});

If $r'=3$, then $\left|X\right|=6$ and $G[X]$ is connected since the minimum degree of $G[X]$ is equal to $\frac{|X|}{2}$. By \cite{regular}, there are two non-isomorphic connected $3$-regular graphs with $6$ vertices (see Figure \ref{agraph}). Since $A_2$ contains $S_1$ as an induced subgraph, by Claim 1, $G[X]\ncong A_2$. Then $G[X]\cong A_1$, and the graph $G_2$ in Figure \ref{mugraph} is the only non-isomorphic graph satisfying (\ref{u1rho}).
\begin{figure}[htb!]
	\centering
		\begin{minipage}{230pt}
	\subfloat[$S_1$]{
		\begin{minipage}{60pt}
			\centering
			\includegraphics[scale=0.1]{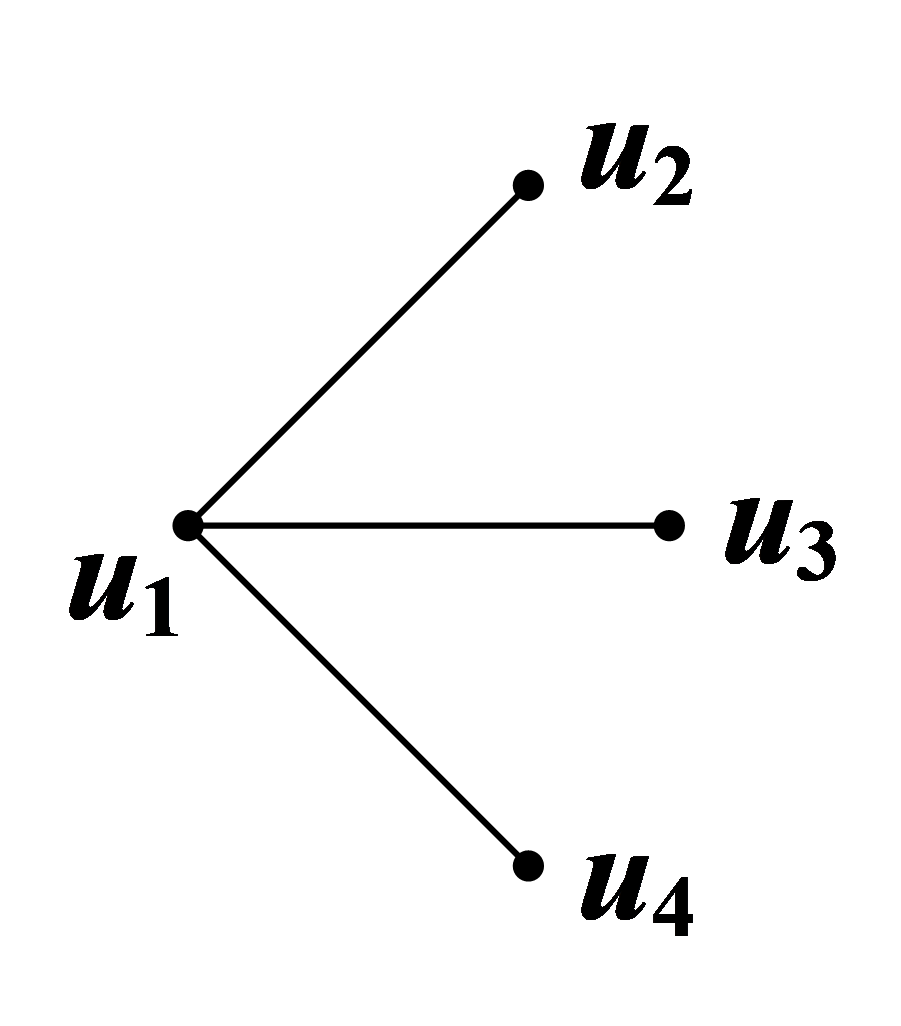}
	\end{minipage}}
	\quad
	\subfloat[$S_2$]{
		\begin{minipage}{65pt}
			\centering
			\includegraphics[scale=0.09]{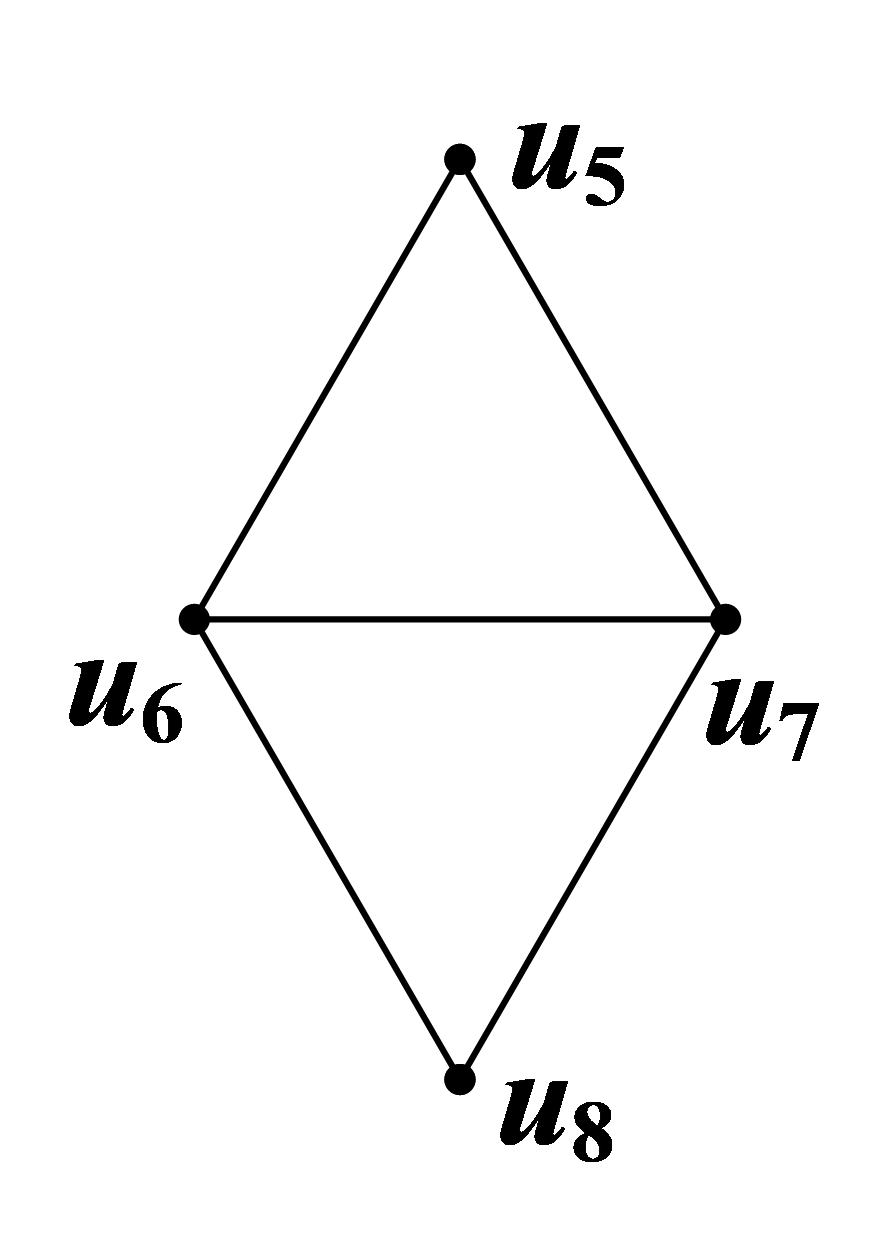}
	\end{minipage}}
	\subfloat[$S_3$]{
		\begin{minipage}{60pt}
			\centering
			\includegraphics[scale=0.1]{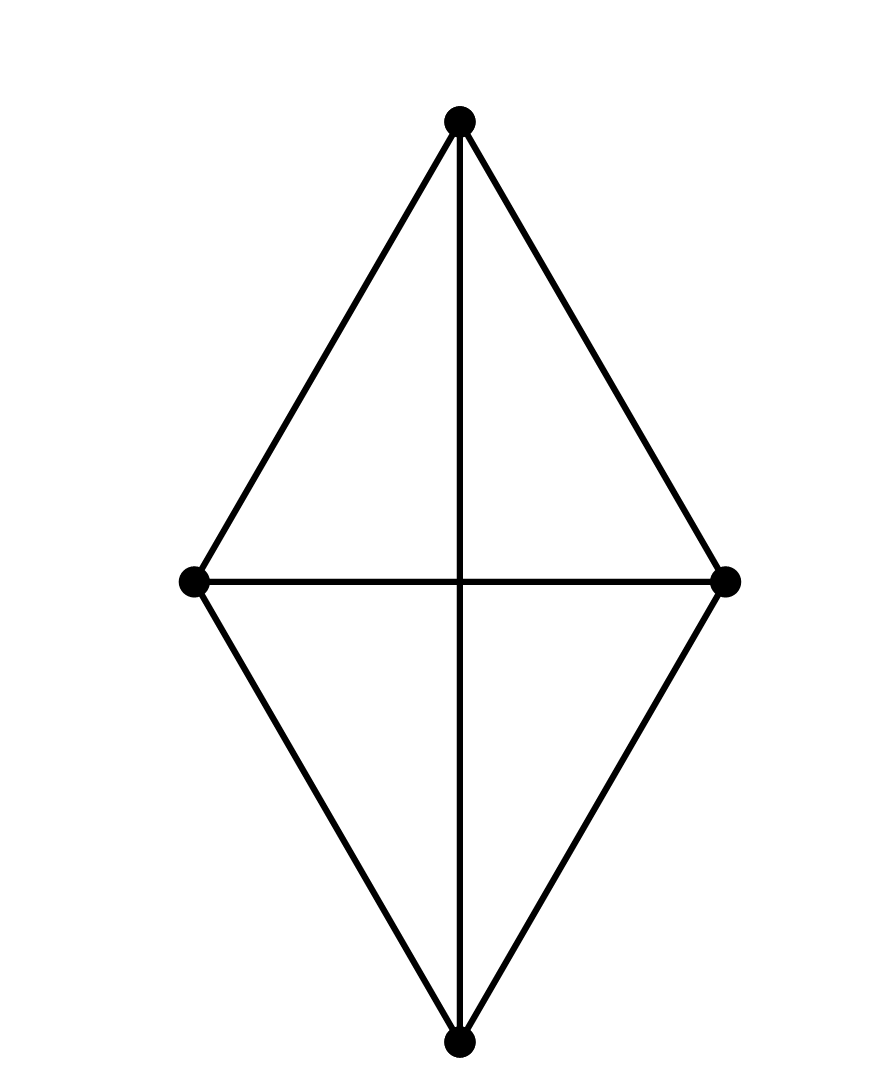}
	\end{minipage}}
	\caption{induced subgraph}
	\label{subgra}
		\end{minipage}
	\begin{minipage}{200pt}
	\centering	
	\subfloat[$A_1$]{
		\begin{minipage}{90pt}
			\centering
			\includegraphics[scale=0.3]{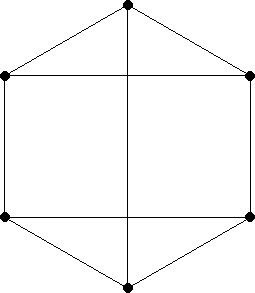}
	\end{minipage}}
	\subfloat[$A_2$]{
		\begin{minipage}{90pt}
			\centering
			\includegraphics[scale=0.3]{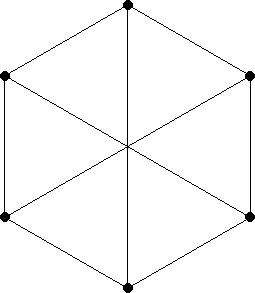}
	\end{minipage}}
	\caption{$3$-regular graphs on $6$ vertices.}
	\label{agraph}
	\end{minipage}
\end{figure}


If $r'=4$, then $\left|X\right|=9$ and $G[X]$ is connected since the minimum degree of $G[X]$ is equal to $\left \lfloor \frac{|X|}{2} \right \rfloor$. By \cite{regular}, there are $16$ non-isomorphic connected $4$-regular graphs with $9$ vertices (see Figure \ref{bgraph}). Since $B_i (i\in \{2,\dots,12\})$ contain $S_1$ as an induced subgraph, graph $B_{13}$ contain $S_3$ as an induced subgraph, graph $B_j (j\in \{14,15,16\})$ contain $S_2$ as an induced subgraph, by Claim 1, we have $G[X]\ncong B_i (i\in \{2,\dots,16\})$. So $G[X]\cong B_1$, and the graph $G_3$ in Figure \ref{mugraph} is the only non-isomorphic graph satisfying (\ref{u1rho}).

Combining the proof of Subcase 2.3, (4) holds.
\begin{figure}[h!]
	\centering
	\subfloat[$B_1$]{
		\begin{minipage}{52pt}
			\centering
			\includegraphics[scale=0.3]{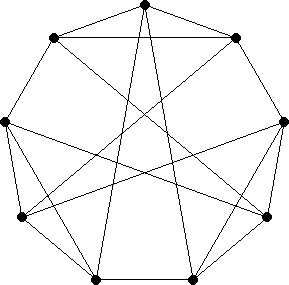}
	\end{minipage}}
	\qquad
	\subfloat[$B_2$]{
		\begin{minipage}{52pt}
			\centering
			\includegraphics[scale=0.3]{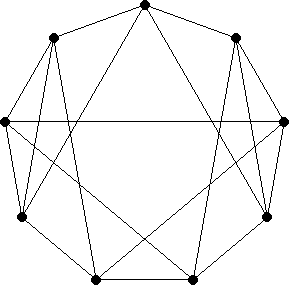}
	\end{minipage}}
	\qquad
	\subfloat[$B_3$]{
		\begin{minipage}{52pt}
			\centering
			\includegraphics[scale=0.3]{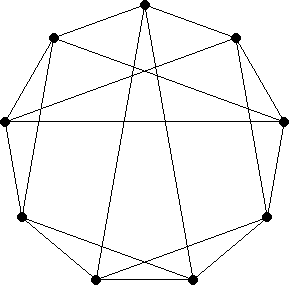}
	\end{minipage}}
	\qquad
	\subfloat[$B_4$]{
		\begin{minipage}{52pt}
			\centering
			\includegraphics[scale=0.3]{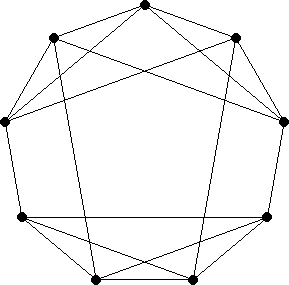}
	\end{minipage}}
	\qquad
	\subfloat[$B_5$]{
		\begin{minipage}{52pt}
			\centering
			\includegraphics[scale=0.3]{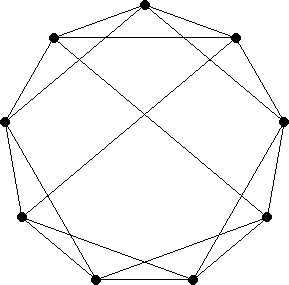}
	\end{minipage}}
	\qquad
	\subfloat[$B_6$]{
		\begin{minipage}{52pt}
			\centering
			\includegraphics[scale=0.3]{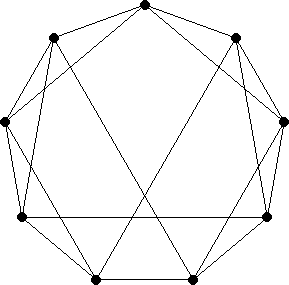}
	\end{minipage}}
	\qquad
	\subfloat[$B_7$]{
		\begin{minipage}{52pt}
			\centering
			\includegraphics[scale=0.3]{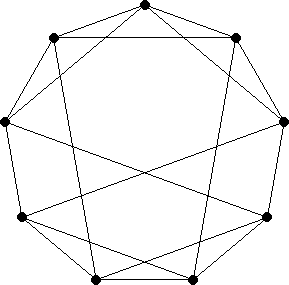}
	\end{minipage}}
	\qquad
	\subfloat[$B_8$]{
		\begin{minipage}{52pt}
			\centering
			\includegraphics[scale=0.3]{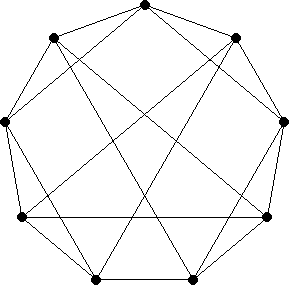}
	\end{minipage}}
	\qquad
	\subfloat[$B_9$]{
		\begin{minipage}{52pt}
			\centering
			\includegraphics[scale=0.3]{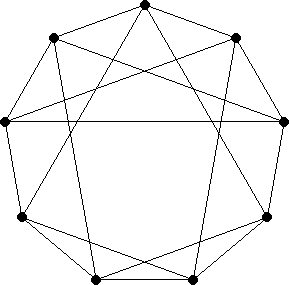}
	\end{minipage}}
	\qquad
	\subfloat[$B_{10}$]{
		\begin{minipage}{52pt}
			\centering
			\includegraphics[scale=0.3]{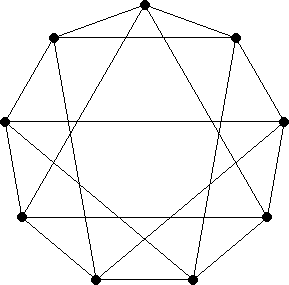}
	\end{minipage}}
	\qquad
	\subfloat[$B_{11}$]{
		\begin{minipage}{52pt}
			\centering
			\includegraphics[scale=0.3]{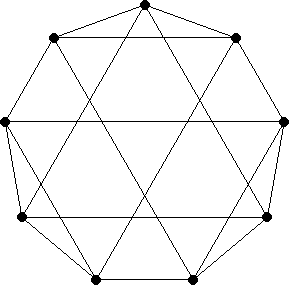}
	\end{minipage}}
	\qquad
	\subfloat[$B_{12}$]{
		\begin{minipage}{52pt}
			\centering
			\includegraphics[scale=0.3]{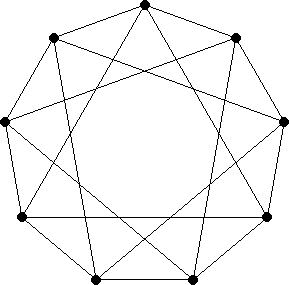}
	\end{minipage}}
	\qquad
	\subfloat[$B_{13}$]{
		\begin{minipage}{52pt}
			\centering
			\includegraphics[scale=0.3]{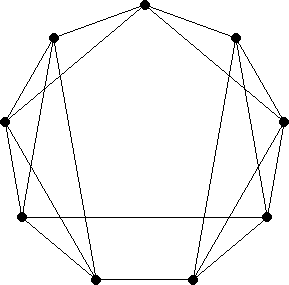}
	\end{minipage}}
	\qquad
	\subfloat[$B_{14}$]{
		\begin{minipage}{52pt}
			\centering
			\includegraphics[scale=0.3]{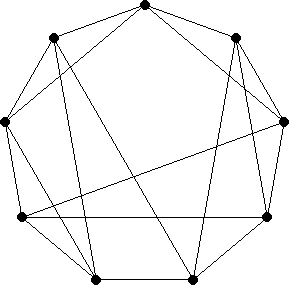}
	\end{minipage}}
	\qquad
	\subfloat[$B_{15}$]{
		\begin{minipage}{52pt}
			\centering
			\includegraphics[scale=0.3]{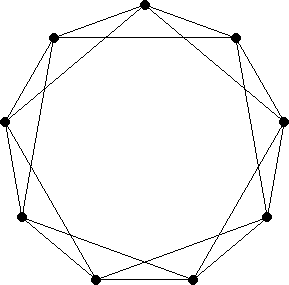}
	\end{minipage}}
	\qquad
	\subfloat[$B_{16}$]{
		\begin{minipage}{52pt}
			\centering
			\includegraphics[scale=0.3]{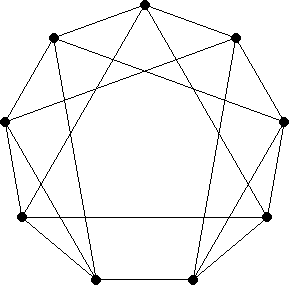}
	\end{minipage}}
	\caption{$4$-regular graphs on $9$ vertices}
	\label{bgraph}
\end{figure}


\noindent\textbf{Case 3.} The vertices in $X$ are of Type III. 

Then $H=K_{3,(\mu^4+7\mu^3+13\mu^2+2\mu-6)/(\mu+2)}$ and all vertices in $X$ are of type $(2,(\mu^3+3\mu^2-2)/(\mu+2))$. Since $\mu \neq -3$, by (\ref{eq27}), we have
$$ \rho_{uv}=\left\{
    \begin{matrix}
    \frac{2}{\mu+2}, & u\sim v,\\
    \frac{\mu^2+2\mu+2}{\mu+2}, & u\nsim v.
    \end{matrix}
    \right.$$
Then $(\mu+2)\mid 2$ by $\rho_{uv}\in  \mathbb{Z}$. Noting that $\rho_{uv}\geq 0$ and $\mu \notin \left\{-1,-3,0\right\}$,
it is impossible. Therefore, there is not an $r$-regular graph with $K_{3,s}$ as a star complement in this case.

Combining the above argument, we complete the proof.$\quad\square$
\end{Proof}

\section{Regular graphs with $K_{s,s}$ as a star complement}\label{kss}
\hspace{1.5em}
By (4) of Theorem \ref{thm35}, we know that there are three regular graphs with $K_{3,3}$ as a star complement for $\mu=1$. In the following, we will study some properties of regular graphs with $K_{s,s}$ as a star complement for an eigenvalue $\mu$, and then give a sharp upper bound for the multiplicity $k$ of $\mu$. 

Since $\mu$ is not an eigenvalue of $K_{s,s}$, we have $\mu\neq 0$ and $\mu \neq \pm s$. When $\mu=-1$, we have $r\equiv -1( {\rm mod}\ (s+1))$ and $G\cong G(r)$ by Theorem \ref{u=-1}. So we discuss the case when $\mu\notin \{-1,0\}$ in the following. 
\begin{prop}
	Let $s\geq 2$ and $G$ be an $r$-regular graph with $K_{s,s}$ as a star complement for an eigenvalue $\mu$, where $\mu\notin \{-1,0\} $. Then
	
	(1) $\mu \in \mathbb{Z}$ and $|\mu|<s$. 
	
	(2) If $\mu=1$, then $s=3$ and $G\cong G_1$, $G_2$ or $G_3$. 
	
	(3) If all vertices in the star set $X$ are of the same type, then $G$ is an r-regular graph of order $(2\mu+1)(\frac{r}{\mu}-1)$. 
\end{prop}
\begin{Proof}
	Clearly, there is no regular graph $G$ with $K_{s,s}$ as a star complement for $\mu=r$ by Theorem \ref{u=r}. Thus $\mu$ is a non-main eigenvalue. Let $t=s$ in (\ref{eq25}), we have $0<a+b=s-\mu\leq2s$, thus $\mu \in \mathbb{Z}$ and $-s\leq \mu <s$. Since $\mu \neq \pm s$, we have $|\mu|<s$, (1) holds.
	
	Since $t=s$, by (\ref{eq25}) and (\ref{eq26}), we have $a=x_1$, $b=x_2$ or $a=x_2$, $b=x_1$ where 
	\begin{equation}\label{solu}
	x_1=\frac{s-\mu +\sqrt{-(s+\mu)(2\mu^2+\mu-s)}}{2},\ x_2=\frac{s-\mu -\sqrt{-(s+\mu)(2\mu^2+\mu-s)}}{2}.
	\end{equation}
Since $a,b\in \mathbb{Z}$, $-(s+\mu)(2\mu^2+\mu-s)=(s-\mu^2)^2-(\mu^2+\mu)^2$ must be a perfect square. Thus, when $\mu=1$, $(s-1)^2-4$ must be a perfect square, so $s=3$, and thus the graphs $G_1$, $G_2$ and $G_3$ (see Figure \ref{mugraph}) are the regular graphs with $K_{3,3}$ as a star complement for $\mu=1$ by (4) of Theorem \ref{thm35}, (2) holds.

Let $u\in X$ be a vertex of type $(a,b)$. Since $G$ is $r$-regular with $K_{s,s}$ as a star complement and all vertices in $X$ are of the same type, we have $a=b$. By (\ref{eq25}) and (\ref{eq26}), we have
$a=b=s=-\mu$ or $a=b=\mu^2, s=2\mu^2+\mu$. 

If $a=b=s=-\mu$, then $|X|=1$, $r=2s=1+s$. Thus $s=1$ and $\mu=-1$, a contradiction. 

If $a=b=\mu^2, s=2\mu^2+\mu$, counting the edges between $X$ and $V(H)$ in two ways, we have $(a+b)\cdot |X|=2s(r-s)$. Thus $|V(G)|=|X|+2s=(2\mu+1)(\frac{r}{\mu}-1)$, (3) holds.$\quad \square$
\end{Proof}

\begin{remark}
	In fact, there are a lot of regular graphs with $H=K_{s,s}$ as a star complement. For example, when $\mu=-2$, it follows from (\ref{solu}) that $(s-4)^2-4$ must be a perfect square, thus $s=2$ or $6$. 
	
	When $s=2$, we have $a=b=2$ by (\ref{solu}). Thus $|X|=1$ by Lemma \ref{lem23}. In this situation, $G$ is not regular. 
	
	But when $s=6$, we have $a=b=4$ by (\ref{solu}). From (\ref{eq27}), we have $\rho_{uv}=\left\{
	\begin{matrix}
	4, & u\nsim v,\\
	6, & u\sim v.
	\end{matrix}
	\right.$ Since $|X|\cdot 8=12(r-6)$, we have $|X|=\frac{3(r-6)}{2} \in \mathbb{Z}$, and then $r$ is even.
	
	Let $(V,W)$ be the bipartition of $H=K_{6,6}$ with $V=\left\{ v_1,v_2,\cdots,v_6 \right\}$,  $W=\left\{ w_1,w_2,\cdots,w_6\right\}$.
	
	If $r=8$, then $|X|=3$. The graph $G_4$ defined as follows is a regular graph with $K_{6,6}$ as a star complement for $\mu=-2$: $X=\{u_1, u_2, u_3\}$, $V(G_4)=V(H)\cup X$, $G_4[X]=3K_1$,
	and $N_H(u_{1})=\left\{v_1,v_2,v_3,v_4,w_1,w_2,w_3,w_4\right\}$,
	$N_H(u_{2})=\left\{v_3,v_4,v_5,v_6,w_3,w_4,w_5,w_6\right\}$,
	$N_H(u_{3})=\left\{v_1,v_2,v_5,v_6,w_1,w_2,w_5,w_6\right\}$. The spectrum of $G_4$ is $[-6,-2^3,0^8,2^2,8]$. 
	
	If $r=10$, then $|X|=6$. The graph $G_5$ defined as follows is a regular graph with $K_{6,6}$ as a star complement for $\mu=-2$: $X=\{u_1, u_2, \dots, u_6\}$, $V(G_5)=V(H)\cup X$, $G_5[X]=C_6$ with the edge set $\{u_1u_2,u_2u_3,u_3u_4,u_4u_5,u_5u_6,u_6u_1\}$,
	and $N_H(u_{1})=\left\{v_1,v_2,v_3,v_4,w_1,w_2,w_3,w_4\right\}$,
	$N_H(u_{2})=\left\{v_2,v_3,v_4,v_5,w_2,w_3,w_4,w_5\right\}$,
	$N_H(u_{3})=\left\{v_3,v_4,v_5,v_6,w_3,w_4,w_5,w_6\right\}$,
	$N_H(u_{4})=\{v_1,v_4,v_5,v_6,w_1, \\ w_4,w_5,w_6\}$,
	$N_H(u_{5})=\left\{v_1,v_2,v_5,v_6,w_1,w_2,w_5,w_6\right\}$,
	$N_H(u_{6})=\left\{v_1,v_2,v_3,v_6,w_1,w_2,w_3,w_6\right\}$. The spectrum of $G_5$ is $[-6,-2^6,0^6,1^2,3^2,10]$. 
	
	Since $|X|\leq \frac{1}{2}(q+1)(q-2)=65$, where $q=|V(H)|$ {\rm(\cite{multi})}, there are various choices of $r$. It can be predicted that there are a lot of graphs that satisfy the conditions. Then the commonalities of the regular graphs with $K_{s,s}$ as a star complement seems to be an interesting question worth studying.	
\end{remark}

It is shown in \cite{mulre} that if $G$ is a connected $r$-regular graph of order $n$ with $\mu \notin \left\{-1,0\right\}$ as an eigenvalue of multiplicity $k$ and $r>2$, $q=n-k$, then $k\leq \frac{1}{2}(r-1)q$. In the following, we will show that when $G$ has $K_{s,s}\ (q=2s\geq 2)$ as a star complement for $\mu$, then $k\leq s(r-s)= \frac{1}{2}q(r-\frac{q}{2})\leq \frac{1}{2}(r-1)q$.


For subsets $U',V'$ of $V(G)$, we write $E(V',U')$ for the set of edges between $U'$ and $V'$. The authors of \cite{match} have determined all the graphs with a star set $X$ for which $E(X,\overline{X})$ is a perfect matching. The result is as follows.
\begin{theo}{\rm (\cite{match})}\label{match}
	Let $G$ be a graph with $X$ as a star set for an eigenvalue $\mu$. If $E(X,\overline{X})$ is a perfect matching, then one of the following holds:\\	{\rm(1)} $G=K_2$ and $\mu=\pm 1$; {\rm(2)} $G=C_4$ and $\mu=0$; {\rm(3)} G is the Petersen graph and $\mu=1$.
\end{theo}

\begin{theo}
	Let $G$ be an $r$-regular graph of order $n$ with $K_{s,s}$ as a star complement for the eigenvalue $\mu \notin \{-1,0\}$ of multiplicity $k$. Then $k\leq s(r-s)$, equivalently $n\leq s(r-s+2)$, with equality if and only if $\mu=1$, $G\cong G_1,G_2$ or $G_3$ {\rm(}see Figure \ref{mugraph}{\rm)}.
\end{theo}
\begin{Proof}
	Let $H=K_{s,s}$ and $V(H)=V\cup W$ with $|V|=|W|=s$. By Lemma \ref{lem23}, $V(K_{s,s})$ is a location-dominating set, and so $G$ is connected. By Lemma \ref{lem23}, we have $|N_H(u)|\geq 1$. Thus $k=|X|\leq \sum_{u\in X} |N_H(u)|=|E(X,\overline{X})|=2s(r-s)$. 
	
	When the equality holds, we have $|N_H(u)|=1$ for all $u\in X$. Since the neighbourhoods $N_H(u)\ (u\in X)$ are distinct, any vertex in $H$ has at most one adjacent vertex in $X$, thus $2s=|\overline{X}|\geq |X|=2s(r-s)$ which means $r\leq s+1$. On the other hand, we have $r\geq s+1$ by $G$ is $r$-regular and connected. Thus $r=s+1$ and then $|X|=2s$. Therefore $E(X,\overline{X})$ is a perfect matching, but there is no such graph $G$ by Theorem \ref{match}.
	
	Let $t=s$ in (\ref{eq25}), we find that $|N_H(u)|=a+b=s-\mu$ is a constant, which means $|N_H(u_1)|=|N_H(u_2)|$ for any $u_1,u_2 \in X$. Therefore, we have $|N_H(u)|\geq 2$ for any $u\in X$ and $2k\leq \sum\limits_{u\in X} |N_H(u)|=|E(X,\overline{X})|=2s(r-s)$, equivalently $n\leq s(r-s+2)$, with equality if and only if $|N_H(u)|=2$ for any $u\in X$.
	
	If $n=s(r-s+2)$, we have $\mu=s-2$ and the possible types for the vertices in $X$ are $(1,1),(0,2),(2,0)$. 
	
	If there are two vertices of type $(0,2)$ (or $(2,0)$) in $X$, then it follows from (\ref{eq27}) that $\rho_{uv}=\left\{
	\begin{matrix}
	\frac{s}{s-1}, & u\nsim v,\\
	\frac{s^2-4s+2}{1-s}, & u\sim v.
	\end{matrix}
	\right.$ By (2) of Lemma \ref{lem23}, we have $\rho_{uv}\neq 2$, then $\rho_{uv}=0$ or $1$, it is a contradiction with $s\in \mathbb{Z}_{+}$. Therefore, there is at most one vertex of type $(0,2)$ (or $(2,0)$).
	
	If there is one vertex of type $(2,0)$ and one vertex of type $(0,2)$ in $X$, then it follows from (\ref{eq27}) that $\rho_{uv}=\left\{
	\begin{matrix}
	\frac{s-2}{s-1}, & u\nsim v,\\
	\frac{s^2-4s+4}{1-s}, & u\sim v.
	\end{matrix}
	\right.$ Clearly, $\rho_{uv}=0$, $s=2$, and $\mu=0$, it is a contradiction. Therefore, the vertex of type $(2,0)$ and the vertex of type $(0,2)$ cannot exist at the same time in $X$.
	
	Suppose that $X$ contains a vertex of type $(0,2)$ (or $(2,0)$), then $|E(X,V)|\neq |E(X,W)|$, a contradiction. Thus all vertices in $X$ are of type $(1,1)$. From (\ref{eq27}), we have $\rho_{uv}=\left\{
	\begin{matrix}
	1, & u\nsim v,\\
	3-s, & u\sim v.
	\end{matrix}
	\right.$ 
	
	If for any $u,v \in X$, $u\nsim v$, then $r=2$ and $H\cong K_{1,1}$, $G\cong C_3$ and thus $s=1$ and $\mu=s-2=-1$, a contradiction. 
	
	If there exists $u,v \in X$ such that $u\sim v$, then $3-s=0$ or $1$ by (2) of Lemma \ref{lem23}. If $3-s=1$, then $s=2$ and $\mu=0$, a contradiction. If $3-s=0$, then $s=3$, $\mu=1$ and thus $G=G_1,G_2$ or $G_3$ {\rm(}see Figure \ref{mugraph}{\rm)} by (4) of Theorem \ref{thm35}. 
	
	The proof is completed. $\quad \square$
\end{Proof}

\section*{Competing interests}

The authors declare that they have no competing interests.

\end{document}